\documentclass[a4paper,12pt]{article}
\usepackage{amsfonts,amssymb,amsmath,amsthm}
\usepackage[dvips]{graphics}
\usepackage[dvips]{graphicx}
\usepackage[dvips]{color}

\def\num{\hspace{-2mm}{\bf }\hspace{2mm}}

\newtheorem{st}{Statement}[section]

\newtheorem{propo}[st]{Proposition}

\newtheorem{thm}[st]{Theorem}
\newtheorem{lemm}[st]{Lemma}

\def\max{{\rm max\,}}
\def\min{{\rm min\,}}

\def\card{{\rm card\,}}

\def\cF{{\cal F}}

\def\cP{{\cal P}}

\def\Fin{{\textbf S}}
\def\card{{\rm card\,}}

\def\Proof:{ \vspace{-1.5mm} {\noindent\it Proof.}}
\def\gm{\vspace{0mm}}
\def\gd{\vspace{0mm}}

\def\Box{\rule{1.5mm}{1.5mm}}

\begin{document}

\renewcommand{\thefootnote}{\fnsymbol{footnote}}

 \title{\gd \gd Thin ultrafilters, P-hierarchy and Martin Axiom}
 \author{ Micha\l \ Machura and Andrzej Starosolski}
 \date{\today}
\maketitle

\gd
\begin{abstract}
Under MA we prove that for the ideal $\cal I$ of thin sets on $\omega$ and for any ordinal $\gamma \leq \omega_1$ there is an ${\cal I}$-ultrafilter (in the sense of Baumgartner), which belongs to the class ${\cal P}_{\gamma}$ of P-hierarchy of ultrafilters. Since the class of ${\cal P}_2$ ultrafilters coincides with a class of P-points, out result generalize theorem of Fla\v{s}kov\'a, which states that there are ${\cal I}$-ultrafilters which are not P-points. It is also related to theorem which states that  
under CH for any tall P-ideal $\cal I$ on $\omega$ there is an ${\cal I}$-ultrafilter, however the ideal of thin sets is not P-ideal. 

\end{abstract}


\gd
\section{Introduction}

Baumgartner in the article \textit{Ultrafilters on $\omega$} (\cite{Baum}) introduced a notion of ${\cal I}$-ultrafilters: 
\vspace{5mm}

A filter on $\omega$ is an ${\cal I}$-ultrafilters, if and only if, for every function $f\in \omega^{\omega}$ there is a set $U\in u $ such that $f[U] \in {\cal I}$.
\vspace{5mm}

 This kind of ultrafilters was studied by large group of mathematician. We shall mention only the most important papers in this subject from our point of view: J. Brendle \cite{Brendle}, C. Laflamme \cite{Laf}, Shelah (\cite{Shelah}). The theory of $\cal I$-ultrafilters was developted by Fla\v{s}kov\'a in a series of articles and in her Ph.D thesis \cite{FlasDok}.  
\vspace{5mm}

In this paper we restrict our attention to \textit{thin ultrafilters} i.e. ${\cal I}$-ultrafilters where  ideal 
${\cal I} $ is ideal  generated by thin sets (\textit{thin ideal}): 

A set $A\subset \omega$ with enumeration $A=\{ a_m : n<\omega  \}$ is called \textit{thin} if $$\lim_{n\to \infty} \frac{a_n}{a_{n+1}} = 0 .$$  
The thin ideal is an example of tall ideal but it is not $P$-ideal.

In \cite{FlasDok} Fla\v{s}kov\'a  proved under Martin Axiom that there are thin ultrafilters, which are not a $P$-point. \vspace{5mm}

Ultrafilters on $\omega$ may be classified with respect to
sequential contours of different ranks, that is, iterations of the
Fr\'{e}chet filter by contour operations. This way an
$\omega_1$-sequence $\{\cP_\alpha \}_{1\leq\alpha\leq\omega_1}$ of
pairwise disjoint classes of ultrafilters - the P-hierarchy - is
obtained, where P-points correspond to the class $\cP_2$, allowing
us to look at the P-hierarchy as the extension of notion of P-point. 
The following theorem was proved by Starosolski, see \cite{Star-P-hier} Proposition 2.1: 

\begin{propo}
An ultrafilter $u$ is a P-point if and only if $u$ belongs to the class ${\cal P}_2$ in P-hierarchy. 
\end{propo}

Many inmportant information about P-hierarchy may be found in \cite{Star-P-hier}.
For additional information regarding sequential cascades and
contours one can look at \cite{DolMyn}, \cite{DolStaWat},
\cite{Dol-multi}, \cite{Star-ff}. However the 
most important definitions and conventions shall be repeated below.

Since $P$-point correspond to ${\cal P}_2$ ultrafilter in P-hierarchy of ultrafilters (more about P-hierarchy one can find below), it would interesting to know to which classes of P-hierarchy can belong ${\cal I}$-ultrafilters. In \cite{Machura-Staros} it was shown that under CH in every class ${\cal P}_{\alpha}$ there are ${\cal I}$-ultrafilters. In this paper we improve this result replacing CH by MA for thin ultrafilters. Let us introduce all necessary definitions and tools. 
\vspace{5mm}

The set of natural numbers (finite ordinal numbers) we denote $\omega$. The filter considered in this paper will be defined on infinite countable set (except one indicated case). This will be usually a set $\max V$ of maximal elements of a cascade $V$ (see definition of cascade below) and we will often identify it with $\omega$ without indication.  The following convention we be applied without mentioning it: 
\vspace{5mm}

If $u$ is a filter on $A \subset B$, then we identify $u$ with the
filter on $B$ for which $u$ is a filter-base. If $\cal F$ is a filter base, then by $\langle \cal F \rangle$ we denote a filter generated by $\cal F$.
\vspace{5mm}

The {\it cascade} is a tree $V$ without infinite branches and with a
least element $\emptyset _V$. A cascade is $\it sequential$ if for
each non-maximal element of $V$ ($v \in V \setminus \max V$) the set
$v^{+V}$ of immediate successors of $v$ (in $V$) is countably
infinite. We write $v^+$ instead of $v^{+W}$ if it is known in which
cascade the successors of $v$ are considered. If $v \in V \setminus
\max V$, then the set $v^+$ (as infinite) may be endowed with an
order of the type $\omega$, and then by $(v_n)_{n \in \omega}$ we
denote the sequence of elements of $v^+$, and by $v_{nW}$ - the
$n$-th element of $v^{+W}$.
\vspace{5mm}

The {\it rank} of $v \in V$ ($r_V(v)$ or $r(v)$) is defined
inductively as follows: $r(v)=0$ if $v \in \max V$, and otherwise
$r(v)$ is the least ordinal greater than the ranks of all immediate
successors of $v$. The rank $r(V)$ of the cascade $V$ is, by
definition, the rank of $\emptyset_V$. If it is possible to order
all sets $v^+$ (for $v \in V \setminus \max V$)  so that for each $v
\in V \setminus \max V$ the sequence $(r(v_n)_{n<\omega})$ is
non-decreasing, then the cascade $V$ is {\it monotone}, and we fix
such an order on $V$ without indication.
\vspace{5mm}

For $v \in V$ we denote by $v^\uparrow$ or $v^{\uparrow V}$ a subcascade of
$V$ built by $v$ and all successors of $v$. 
\vspace{5mm}

If $\mathbb{F}= \{{\cal F}_s: s \in S \}$ is a family of filters on
$X$ and if $\cal G$ is a filter on $S$, then the {\it contour of $\{
{\cal F}_s \}$ along $\cal G$} is defined by
$$\int_{\cal G} \mathbb{F} = \int_{{\cal G}} \{ {\cal F}_s  : s \in S \}  =
\bigcup_{G \in {\cal G}} \bigcap_{s \in G} {\cal F}_s.$$
\vspace{5mm}

Such a construction has been used by many authors (\cite{Fro},
\cite{Gri1}, \cite{Gri2}) and is also known as a sum (or as a limit)
of filters.

Operation of sum of filters we apply in definition of \textit{contour of cascade}:
Fix a  cascade $V$. 
Let ${\cal G}(v)$ be a filter on $v^{+}$ for every $v \in V \setminus \max V$.
For $v \in \max V$ let ${\cal G}(v)$ be a trivial ultrafilter on a singleton $\{ v \}$ (we can treat it as principal ultrafilter on $\max V$ according to convention we assumed).
This way we have defined a function $v \longmapsto  {\cal G}(v)$. 
 Contour of each sub-cascade $v^{\uparrow}   $ is defined inductively with respect to rank of $v$:
$$\int^{\cal G} v^{\uparrow} = \{ \{ v \} \}$$ 
 for $v\in \max V$  (i.e. $\int^{\cal G} v^{\uparrow}$ is just  a trivial ultrafilter on singleton $\{ v \} $) ; 
$$\int^{\cal G} v^{\uparrow} =\int_{{\cal G}(v)} \left\{ \int^{\cal G} w^{\uparrow} : w\in v^{+} \right\} $$ for $v \in V \setminus \max v$.
Eventually we put 
$$ \int^{\cal G} V = \int^{\cal G} \emptyset_{V}.$$
Usually we shall assume that all the filters ${\cal G}(v)$ are Frechet (for $v\in V\setminus \max V$). In that case we shall write $\int V$ instead of $\int^{\cal G} V$.

Similar
filters were considered in \cite{Kat1}, \cite{Kat2}, \cite{Dagu}.
Let $V$ be a monotone sequential cascade and let $u=\int V$.
Then a {\it rank $r(u)$} of $u$ is, by definition, the rank of $V$.

It
was shown in \cite{DolStaWat} that if $\int V= \int W$, then $r(V) =
r(W)$.
\vspace{5mm}

We shall say that a set $F$ \textit{meshes} a contour ${\cal V}$ ($F \# {\cal V}$) if and only if ${\cal V} \cup \{ F\}$ has finite intersection property and can be extended to a filter. If $\omega \setminus F \in {\cal V}$, then we say that $F$ is \textit{residual} with respect to ${\cal V}$ .
\vspace{5mm}

Let us define ${\cP}_\alpha$
for $1 \leq {\alpha <\omega_1}$ on $\beta\omega$ (see
\cite{Star-P-hier}) as follows: $u \in {\cP}_\alpha$ if there is no
monotone sequential contour $C_{\alpha}$ of rank $\alpha$ such that
$C_{\alpha}  \subset u$, and for each $\beta$ in the range
$1\leq\beta < \alpha$ there exists a monotone sequential contour
$C_\beta$ of rank $\beta$ such that $C_\beta \subset u$. Moreover,
if for each $\alpha < \omega_1$ there exists a monotone sequential
contour $C_\alpha$ of rank $\alpha$ such that $C_\alpha \subset u$,
then we write $u \in {\cP}_{\omega_1}$.
\vspace{5mm}

Let us consider a monotone cascade $V$ and a monotone sequential
cascade $W$. We will say that $W$ is a sequential extension of $V$
if:

1) $V$ is a subcascade of cascade $W$,

2) if $v^{+V}$ is infinite, then $v^{+V} = v^{+W}$,

3) $r_V(v) = r_W(v)$ for each $v \in V$.

Obviously, a monotone cascade may  have many sequential extensions.

Notice that if $W$ is a sequential extension of $V$ and $U \subset
\max V $, then $U$ is residual for $V$ if and only if $U$ is
residual for $W$.
\vspace{5mm}

The following theorem was proved in \cite{Star-P-hier} Theorem 2.8:
\begin{thm}\num The following statements are equivalent:
\begin{enumerate}
\item P-points exist,
\item $\cP_\alpha$ classes are non-empty for each countable successor
$\alpha$,
\item There exists a countable  successor $\alpha>1$ such that the
class $\cP_\alpha$ is non-empty.
\end{enumerate}
\end{thm}

Starosolski has proved in \cite{Star-Top-P-Hier} Theorem 2.6  that under CH every class ${\cal P}_{\alpha}$ is nonempty. 
\vspace{5mm}


\section{Lemmas}

The following lemmas will be used in the prove of a main theorem. 

The first lemma is one of  lemmas proved in \cite{Star-Top-v-P-Hier} (see: Lemma 6.3 ):

 \begin{lemm}\num
Let $\alpha < \omega_1$ be a limit ordinal and let $({\cal V}_n : n<\omega)$ be
a sequence of monotone
sequential contours such that $r({\cal V}_n)<r({\cal
V}_{n+1})<\alpha$ for every $n$ and that $\bigcup_{n < \omega}{\cal V}_n$ has finite intersection property. 
Then there is no monotone sequential contour
$\cal{W}$ of rank $\alpha$ such that ${\cal{W}} \subset
\langle \bigcup_{n<\omega} {\cal V}_n \rangle$.
\end{lemm}

As a corollary we get: 

\begin{lemm} Assume Martin Axiom. 
Let $\alpha < \omega_1$ be a limit ordinal, let $({\cal V}_n)_{n<\omega}$ be an increasing ("$\subset$")
sequence of monotone sequential
contours, such that $r({\cal V}_n)<\alpha$ and let ${\cal F}$ be a family of
sets of cardinality less than ${\frak c}$ such that $\bigcup_{n<\omega} {\cal V}_n  \cup {\cal F}$ has finite intersection property. Then
$\langle \bigcup_{n<\omega}{\cal V}_n  \cup {\cal F}) \rangle$ does not contain any monotone sequential contour of rank $\alpha$.
\end{lemm}

\textit{Proof:} 
One may assume that ${\cal F}$ is closed under finite intersection and contains Frechet filter.
We consider partially ordered set:
$$\mathbb{P} = \{ \langle S,N \rangle   : S\subset \omega \mbox{ is finite, } N\in {\cal F}  \mbox{ and } \max S < \min N \}.  $$
To simplyfize notion for $p =\langle S,N \rangle \in \mathbb{P}$ we put $S_p=S$ and $N_p= N$. 
An order on $\mathbb{P}$ is given by formula: 
$$ q\leq p \Leftrightarrow S_q \supseteq S_p, \ N_q \subseteq N_p \mbox{ and } S_q\setminus S_p \subset N_p.$$
For a given cascade $V$, finite set $S$ and natural number $k$ we define property $\Phi (V,S,k)$. The definition is inductive: for any $v\in \max V$ we put
$$ \Phi (v,V,S,k) \Leftrightarrow v\in S.$$
Next we define $\Phi (v,V,S,k)$ for  $v\in V\setminus \max V$: 
$$ \Phi (v,V,S,k) \Leftrightarrow \card \{ w\in v^{+} : \Phi (w,V,S,k)  \} \geq k.$$
Eventually we put $ \Phi (V,S,k) = \Phi (\emptyset_V,V,S,k)$.
\vspace{0.5cm}
Let $V_n$ be cascades such that $\int V_n = {\cal V}_n$, for any $n$. 
We define dense sets for every ${\cal V}_n$ and $k<\omega$: 
$$\mathbb{D}_{n,k} =\{ p \in \mathbb{P} :   \Phi (V_n,S_p,k)   \}; $$
and for every $F\in {\cal F}$: 
$$\mathbb{D}_{F} =\{ p \in \mathbb{P} :   N_p \subset F   \}; $$ 
Let $\mathbb{G}$ be a generic filter and let $G= \bigcup \{ S_p : p\in {\mathbb G} \}$.
It is easy to see that $G$ is infinite, $G \setminus F$ is finite for every $F\in {\cal F}$ and 
$\bigcup_{n<\omega} {\cal V}_n  \cup {\cal F} \cup \{ G\} $ has finite intersection property.
\vspace{0.5cm}

Let ${\cal W}_n = \{ U \cap G : U\in {\cal V}_n \}$ for every $n$. It is easy to see that ${\cal W}_n$ is monotone sequential contour of the same rank as ${\cal V}_n$. Consider a sequence $({\cal W}_n :n<\omega)$. By Lemma 2.1 the union $\bigcup_{n<\omega} {\cal W}_n$ do not contains contour of rank $\alpha$.

 $\hfill$ $\Box$
\vspace{5mm}

Similarly but easer one can prove:

\begin{lemm} Assume Martin Axiom. 
Let $\alpha = \delta +1 $ be a succesor ordinal, let ${\cal V}$ be 
sequence of monotone sequential
contour of rank $r({\cal V})= \delta$ and let ${\cal F}$ be a family of
sets of cardinality less than ${\frak c}$ such that $ {\cal V}  \cup {\cal F}$ has finite intersection property. Then
$\langle {\cal V}  \cup {\cal F}) \rangle$ does not contain any monotone sequential contour of rank $\alpha$.
\end{lemm}

\section{Main result}

In this section we shall present main result of the paper. 

\begin{thm} Assume Martin Axiom. Let $\cal I$ be thin ideal, 
and let $0<\gamma <\omega_1$ be an ordinal.
Then there exists an $\cal I$-ultrafilter $u$ which belongs to $\cP_\gamma$.
\end{thm}
\Proof:
We shall split proof into four cases: $\gamma=1$, $\gamma=2$, $\gamma >2$ is a succesor ordinal (the main step), $\gamma$ is limit ordinal. 
\vspace{0.5cm}

\textbf{Step 0:} $\gamma =1$ -trivial.
\vspace{0.5cm}

\textbf{Step 1:} $\gamma =2$. Let $\{{\cal W}_{\alpha}, \alpha<\omega_1\}$ be an enumeration of all monotone sequential contours
of rank $2$.
Let $\omega^\omega = \{f_\alpha: \alpha< \omega_1\}$.
\hspace{0.5cm}

By transfinite induction, for $\alpha < \omega_1$ we build filter bases  ${\cal F}_\alpha$,
such that:
\begin{enumerate}
\item ${\cal F}_0$ is a Frechet filter;

\item for each $\alpha < \omega_1$ $ card ({\cal F}_\alpha) = \alpha \cdot \omega $;

\item ${\cal F}_\alpha \subset {\cal F}_\beta$ for $\alpha<\beta$;

\item ${\cal F}_\alpha = \bigcup_{\beta<\alpha}{\cal F}_\beta$ for $\alpha$ limit ordinal;


\item for each $\alpha<\omega_1$ there is $ F \in {\cal F}_{\alpha+1}$ such that $f_\alpha[F]\in {\cal I}$;

\item for each $\alpha<\omega_1$ there is $F\in{\cal F}_{\alpha+1}$ such that the complement of $F$ belongs to ${\cal W}_{\alpha}$.
\end{enumerate}

Suppose that ${\cal F}_\alpha$ is already define, we will show how to build ${\cal F}_{\alpha+1}$.
First we shall add to ${\cal F}_{\alpha}$ a set $G_{\alpha}$ which should  
 take care on conditions 5. 
Next we will take care on condition 6 by adding set $A_{\alpha}$ to the list of generators of ${\cal F}_{\alpha +1}$.

If there is $U$ such that $f_{\alpha}[U]$ is finite and ${\cal F}_{\alpha}\cup \{ U \}$ has finite intersection property, then we put $G_{\alpha}=U$ and we are done. 
Assume the opposite. In this case we
 shall use Martin's Axiom. Consider a poset:
$$\mathbb{P} =\{ K \subset \omega : K \mbox{ is finine and }  [u,v] \cap K =\{ u, v \}  \Rightarrow v>u^2  \}.  $$ 
For every $F$ and $k<\omega$ we define a dense sets
$$ \mathbb{D}_{F,k} = \{ K: \card ( K\cap f_{\alpha}[F] ) \geq k \}.$$
There is a generic filter $\mathbb{G}$ which intersect each $\mathbb{D}_{F,k}$
Let $G_{\alpha} = f_{\alpha}^{-1} ( \bigcup \{ K :K \in \mathbb{G} \} )$.
\vspace{0.5cm}

Notice that ${\cal F}_{\alpha} \cup \{G_{\alpha} \}$ has finite intersection property and  cardinality less than ${\frak c}$.
A subbase of any sequential contour of rank 2 has cardinality at least ${\frak d} >\aleph_0$, but ${\frak d}= {\frak c}$ under MA, thus none of them one, in partucular ${\cal W}_{\alpha}$, is contained in ${\cal F}_\alpha \cup \{G_{\alpha} \}$.
This means that there is a set $A_{\alpha}$ such that its complement belongs to  ${\cal W}_{\alpha}$ and a family  ${\cal F}_{\alpha +1} = {\cal F}_{\alpha} \cup \{G_{\alpha}, A_{\alpha} \}$
has finite intersection property. 
\vspace{0.5cm}

\textbf{Step 3:} for limit $\gamma$.
The proof in this case is base on the same idea as step 1, but it is more sophisticated and technical. 

Let $({\cal V}_n)_{n<\omega}$ be an increasing ("$\subset$")
sequence of monotone sequential
contours, such that their ranks $r({\cal V}_n)$ are smaller than $\gamma$ but converging to $\gamma$. For each $n< \omega$ denote by $V_n$ a
(fixed) monotone sequential cascade
such that $\int V_n = {\cal V}_n$.
Let $\{{\cal W}_{\alpha}, \alpha< {\frak c} \}$ be an enumeration of all monotone sequential contours
of rank $\gamma$.
Let $\omega^\omega = \{f_\alpha: \alpha< {\frak c} \}$.
\hspace{0.5cm}

By transfinite induction, for $\alpha < {\frak c} $ we build filter bases  ${\cal F}_\alpha$,
such that:
\begin{enumerate}
\item ${\cal F}_0$ is a Frechet filter;

\item for each $\alpha < {\frak c}$, $ card ({\cal F}_\alpha) = \alpha \cdot \omega $;

\item ${\cal F}_\alpha \subset {\cal F}_\beta$ for $\alpha<\beta$;

\item ${\cal F}_\alpha = \bigcup_{\beta<\alpha}{\cal F}_\beta$ for $\alpha$ limit ordinal;

\item $\bigcup_{i<\omega}{\cal V}_i \cup \bigcup_{\alpha<\omega_1}{\cal F}_\alpha$
has finite intersection property;

\item for each $\alpha<\omega_1$ there is $ F \in {\cal F}_{\alpha+1}$ such that $f_\alpha[F]\in {\cal I}$;

\item for each $\alpha<\omega_1$ there is $F\in{\cal F}_{\alpha+1}$ such that the complement of $F$ belongs to ${\cal W}_{\alpha}$.
\end{enumerate}

Suppose that ${\cal F}_\alpha$ is already define, we will show how to build ${\cal F}_{\alpha+1}$.
This shall be done in 2 substeps. 
First step is to fint   
by one set $G_{\alpha}$ such that 
$\bigcup_{m<\omega} {\cal V}_n \cup {\cal  F}_{\alpha} \cup \{ G_{\alpha} \}$ has finite intersection property and $f_{\alpha}[G_{\alpha}] \in {\cal I}$.  For this purpose we use Martin's Axiom.
The set $G_{\alpha}$  take care on all the contours ${\cal V}_n$. Adding it as generator  to $F_{\alpha +1}$ will ensure preservation of conditions 5 and 6. 
On the last step will take care on condition 7 by adding set $A_{\alpha}$ to the list of generators of $F_{\alpha +1}$.

\vspace{0.5cm}

\noindent \textit{Substep i)}  
Let us introduce an axillary definition.
\vspace{0.5cm}

\noindent
\textit{Definition:} Fix a monotone sequential cascade $V$, a set $F$ and  a function $f \in \omega^{\omega}$. For each $v\in V$, we write $U\in \Fin(v)$ if
\begin{enumerate}
\item  $U \subset \max v^\uparrow$;

\item $(U \cap F) \# \int v^\uparrow$;

\item $\card(f[U \cap F])=1$.
\end{enumerate}

We following two claims are crucial, for proof see \cite{Machura-Staros} Proposition 3.2:

\begin{propo}
One that one of the following possibilities holds:
\vspace{5mm} 
 
A) $\Fin(\emptyset_{V})\not=\emptyset$;
\vspace{5mm}

B) there is an antichain (with respect to the order of a cascade) $\mathbb{A}\subset V$ such that:
\begin{enumerate}
\item  $\Fin(v) = \emptyset$ for all $v \in \mathbb{A}$,

\item $ \left( \bigcup \{  \max w^\uparrow   : w\in v^+, \Fin(w)\not= \emptyset  \} \right)  \# \int v^\uparrow$ for all $v \in \mathbb{A}$,

\item $\left( \bigcup \{
\max v^\uparrow :   v \in {\mathbb{A}}  \} \right) \# \int V$.
\end{enumerate}
\end{propo}
\hfill $\Box$


For a given cascade $V$, base of a filter ${\cal F}$ and a function $f\in \omega^{\omega}$.
Put 
$$ J_F^V \{ v: card (f [\max v^{\uparrow} \cap F ] ) =\omega   \}.$$
If a cascade is fixed we write $J_F$ instead of $J_F^V$.

\begin{propo} Let $V$ be a cascade.
Either for every $F\in \cF$ 
$$\bigcup \{ \max v^{\uparrow} : v\in J_F \} \# \int V $$
either there is $H$ such that $f[H]$ finite and $H \# \int V$.
\end{propo}

\textit{Proof of the proposition:}

Assume that for some $F$ 
$$\neg ( \bigcup \{ \max v^{\uparrow} : v\in J_F \} \# \int V )$$
We use previous proposition lemma . There are two possibilities: 

A) There is $U$ such that $(U \cap F) \#  \int V$ and $f[U\cap F]$ is finite. Then put $H=U\cap F$.

B) There is an antichain $\mathbb{A}$ such that:
\begin{itemize}
\item $\bigcup \{ \max v^{\uparrow} : v\in \mathbb{A} \}  \# \int V$
\item $\bigcup \{ \max w^{\uparrow} : w\in v^{+} \}  \# \int v^{\uparrow}   $ for any $v\in \mathbb{A}$. 
\item For every $ v\in \mathbb{A}  $  there is NO $T_v$ 
\begin{itemize}
\item $T_v \subset \max v^{\uparrow}$
\item $(T_v \cap  F )\# \int v^{\uparrow}$
\item $f[T_v \cap F] $ finite
\end{itemize}
\item For every $ v\in \mathbb{A}  $ and infinite many $w \in v^{+}$ there is $U_w$ 
\begin{itemize}
\item $U_w \subset \max w^{\uparrow}$
\item $(U_w\cap  F )\# \int w^{\uparrow}$
\item $f[U_w \cap F] $ finite
\end{itemize}
\end{itemize}

Put 
$$T_v =\bigcup_{w\in v^{+}} U_w.$$
Observe that 
$f[T_v \cap F]$ has to be infinte. 
But 
$\bigcup_{v\in \mathbb{A}} (T_v \cap  F )\# \int V$
Thus 
$\bigcup_{v \in \mathbb{A}} \max v^{\uparrow}\# \int V$.
On the other hand $\mathbb{A} \subset J_F$, so  $\bigcup_{v \in J_F} \max v^{\uparrow}\# \int V$. 

Contradiction!

\hfill $\Box$

Fix cascades $V_n$ such that $\int V_n= {\cal V}_n$, a function $f=f_{\alpha}$. 
There are two possibilities: 
\vspace{5mm}

\textit{Case A:}
There are $F$ and  infinite many $n$ for which there is $U_n$ such that $(U_n \cap F) \#  \int V_n$ and $f[U_n\cap F]$ is finite. 
Without lost of genereality one may assume that for every $n$ there is such $U_n$ and by shrinking it that $f[U_n\cap F]$ is one-point. 

We put $H_n=U_n\cap F$ and consider two subcases: 
\vspace{0.5cm}

If $\bigcup_n f[H_n]$ is finite then put $G_{\alpha}= f^{-1}[ \bigcup_n f[H_n] ]$. It is obvoius that  
$f[G_{\alpha}]$ is thin as finite and $ \bigcup_{n} {\cal V}_n \cup {\cal F}_{\alpha} \cup \{ G_{\alpha} \} $ has finite intersection property. 
\vspace{0.5cm}

Assume that $\bigcup_n f[H_n]$ is infinite. Let $f[H_n] =\{ x_n  \}$. We can select $x_n$ such that they are far away one from each other and get set $\{ x_n : n \in M\} \in {\cal I}$. Put $G_{\alpha} = f^{-1}[\bigcup_{n\in M} f[H_n]]= f^{-1}[ \{ x_n : n \in M\} ]  $. It is easy to see that $ \bigcup_{n} {\cal V}_n \cup {\cal F}_{\alpha} \cup \{ G_{\alpha} \} $ has finite intersection property, because the union $\bigcup_{n} {\cal V}_n$ is increasing.
\vspace{5mm}
 
\textit{Case B:}
For every $V_n$ and every $F$ we have $$\bigcup \{ \max v^{\uparrow} : v\in J_F^{V_n} \} \# \int V_n $$
In the case B we shall use Martin's Axiom. Consider again the poset:
$$\mathbb{P} =\{ K \subset \omega: K \mbox{ is finite and }   [u,v] \cap K =\{ u, v \}  \Rightarrow v>u^2  \},  $$  
ordered by inclusion. 
\vspace{5mm}

For every $V_n$ and $F$ let $\mathbb{A}^{V_n}_F$ be an antichain defined as in proposition 3.2 . 
For a given $F$ and $v\in   \mathbb{A}^{V_n}_F$ let
$$C^{V_n,F}_{v} = \left\{ w\in v^{+V_{n}} : (\exists U_w )(\exists x_w) \ (U_w\cap F) \# \int v^{\uparrow V_n} \mbox{ and } f(U_w \cap F ) = \{ x_w \} \right\}.$$
Let $X^{V_n,F}_v= \{ x_w : w\in C^{V_n,F}_v \}$ and define a dense sets
$$ \mathbb{D}_{V_n,F,v,k} = \{ K: \card ( K\cap f[\max v^{\uparrow V_n}\cap X^{V_n,F}_v ] ) \geq k \}$$
for every  $k<\omega$.
There is a generic filter $\mathbb{G}$ which intersect each $\mathbb{D}_{V_n,F,v,k}$
Let $G_{\alpha} = f^{-1} ( \bigcup \{ K :K \in \mathbb{G} \} )$.
\vspace{0.5cm}

\noindent \textit{Substep ii)} 
Since the family ${\cal F}_\alpha \cup \{G_\alpha\}$ has cardinality less than ${\frak c}$, thus by Lemma 2.2 
 there exists $A_\alpha$ residual for the contour ${\cal W}_{\alpha}$ and such that a family
 $\bigcup_{n<\omega}{\cal V}_n \cup {\cal F}_\alpha
\cup\{G_\alpha , A_\alpha\}$ has finite intersection property.  Let $ {\cal F}_{\alpha + 1}= {\cal F}_\alpha
\cup\{G_\alpha\} \cup \{A_\alpha\}$  \vspace{0.5cm}

Take any ultrafilter $u$ that extends
$\bigcup_{n<\omega}{\cal V}_n \cup \bigcup_{\alpha<\omega_1}{\cal F}_\alpha$.
By condition 5) $u$ is an $\cal I$-ultrafilter, by condition 6)
$u$ do not contain any monotone sequential contour of rank $\gamma$ which jointly with 
$\bigcup {\cal V}_n \subset u$ give us $u \in {\cal P}_\gamma$.

So the proof is done also for limit $\gamma$.  
\vspace{5mm}

\textbf{Step 2:} $\gamma = \delta +1$ is successor ordinal greater than 2.

The proof is similar as in Step 3, but instead of fixing a sequence of contours $({\cal V}_n : n<\omega )$ it is sufficient to take one contour ${\cal V}$ of rank $\delta$ and cascade $V$ such that $\int V= {\cal V}$. Proof will split into 2 cases: 
\vspace{5mm}

\textit{Case A:} There is $F$ and $U$  such that $(U\cap F) \# \int V$ and $f[U\cap F ]$ is finite. Then put $G_{\alpha}= U\cap F$.

\vspace{5mm}

\textit{Case B:} $\bigcup \{ \max v^{\uparrow} : v\in J^V_F \} \# \int V$.
We proceed exactly like in Case B of Step 3: consider the same poset $\mathbb{P}$ and dense sets $\mathbb{D}_{V,F,v,k}$ to define $G_{\alpha}$

\vspace{5mm}

At the end instead of using Lemma 2.2 one gas to use Lemma 2.3.

\hfill $\Box$

 \hspace{-2mm}
\gd

\medskip

{\small\sc \noindent {Micha\l} Machura, Instytut Matematyki, Uniwersytet {\'{S}l\c{a}ski},
, Katowice, Poland

E-mail:  machura@math.us.edu.pl}

{\small\sc \noindent Andrzej Starosolski, {Wydzia\l}  {Matematyczno-Fizyczny},
Politechnika \'{S}l\c{a}ska, Gliwice, Poland

E-mail:  andrzej.starosolski@polsl.pl}

\end{document}